\documentclass[a4paper,11pt,reqno]{amsart}
\usepackage{amsmath,a4,amsfonts,amssymb,amsthm,enumerate,fleqn,breqn,dsfont}
\usepackage[margin=1in]{geometry}
\usepackage{etoolbox}

\usepackage{enumitem}
\setdescription{leftmargin=0cm,labelindent=0cm}

\def\CC{\mathbb{C}}

\def\RR{\mathbb{R}}
\def\NN{\mathbb{N}}

\def\QQ{\mathfrak{Q}}

\def\T{\mathfrak{T}}
\def\cc#1{\{#1\}}
\def\pp#1{\|#1\|}
\def\ra{\rangle}
\def\la{\langle}
\def\nh{\mathcal{H}}

\def\Q{\mathbb{H}^R(\QQ)}

\def\cc#1{\{#1\}}
\def\pp#1{\|#1\|}
\def\ra{\rangle}
\def\la{\langle}

\theoremstyle{plain}
\newtheorem{theorem}{\bf Theorem}

\theoremstyle{remark}
\newtheorem{definition}[theorem]{\bf Definition}

\newtheorem{example}[theorem]{\bf Example}
\newtheorem{remark}[theorem]{\bf Remark}

\def\Qi{\mathbb{H}^R_i(\QQ)}
\def\xmn{\cc{u_i}_{i\in \NN}}
\def\hmn{\cc{\mathbb{H}^R_i(\QQ)}_{i\in \NN}}
\def\h2n{\cc{\Q_{n,i}}_{i=1,2, \cdots, 2n \atop{n\in \NN}}}
\def\tmn{\cc{\mathfrak{T}_{n,i}}_{i=1,2, \cdots, m_n \atop{n\in \NN}}}
\def\Rmn{\cc{\mathfrak{R}_{n,i}: \Q \to \Q_{n,i}}_{i=1,2, \cdots, m_n \atop{n\in \NN}}}
\def\Tmn{\cc{\mathfrak{T}_{i}: \Q \to \mathbb{H}^R_{ i}(\QQ)}_{i\in \NN}}
\def\ui{\cc{\mathfrak{U}_{i}}_{i\in \NN}}
\def\Umn{\cc{\mathfrak{U}_{i}: \Q \to \Qi}_{i\in \NN}}
\def\Tnn{\cc{\mathfrak{T}_{n,i}: \Q \to \Q_{n,i}}_{i=1,2, \cdots, 2n \atop{n\in \NN}}}
\def\ymn{\cc{v_i}_{i\in\NN}}
\def\xxmn{\cc{\la x, x_{n,i}\ra}_{i=1,2, \cdots, m_n \atop{n\in \NN}}}
\def\alphamn{\cc{\alpha_{n,i}}_{i=1,2, \cdots, m_n \atop{n\in \NN}}}
\def\betamn{\cc{\beta_{n,i}}_{i=1,2, \cdots, m_n \atop{n\in \NN}}}
\def\llim{\lim\limits_{n\to \infty}}
\def\S{{S}}
\def\T{\mathfrak{T}_{\cc{m_n}}}
\def\suml{\sum\limits_{i=1}^{\infty}}
\def\mnsum{\sum\limits_{i=1}^{\infty}}
\def\mpsum{\sum\limits_{i=1}^{m_p}}
\def\mqsum{\sum\limits_{i=1}^{m_q}}
\def\elii{\ell^2(\QQ)}

\title[F\lowercase{rame of Operators in quaternionic} H\lowercase{ilbert spaces}]{Frame of Operators in quaternionic Hilbert spaces}

\author[ S\lowercase{harma}, J\lowercase{arrah, and} K\lowercase{aushik}]
{S.K. S{harma}, A.M. J{arrah}, and S.K. K{aushik}}

\address[S. K. Sharma and S. K. Kaushik]{Department of Mathematics,
Kirori  Mal College,
University of Delhi, Delhi~110~007, INDIA.}
\email[S.K. Sharma]{\lowercase{sumitkumarsharma@gmail.com}}
\email[S.K. Kaushik]{\lowercase{shikk2003@yahoo.co.in}}

\address[A.M. Jarrah]{Department of Mathematics,
	Yarmouk University, Irbid, Jordan.}
\email[]{\lowercase{ajarrah@yu.edu.jo}}


\begin{document}
\maketitle \baselineskip14pt
\begin{abstract}
In this paper,
we introduce and study frame of operators in quaternionic Hilbert spaces as a generalization of $g-$frames which in turn
generalized various notions like Pseduo frames, bounded quasi-projectors and frame of subspaces (fusion frames) in separable  quaternionic Hilbert spaces.
\end{abstract}
\bigskip\bigskip\bigskip

\def\h2n{\cc{\Q_{n,i}}_{i=1,2, \cdots, 2n \atop{n\in \NN}}}
\def\tmn{\cc{\mathfrak{T}_{i}}_{i\in \NN}}
\def\Rmn{\cc{\mathfrak{R}_{i}: \Q \to \Qi}_{i\in \NN}}
\def\rmn{\cc{\mathfrak{R}_{i}}_{i\in \NN}}

\def\Tnn{\cc{\mathfrak{T}_{n,i}: \Q \to \Q_{n,i}}_{i=1,2, \cdots, 2n \atop{n\in \NN}}}
\def\ymn{\cc{y_{n,i}}_{i=1,2, \cdots, m_n \atop{n\in \NN}}}
\def\xxmn{\cc{\la x, x_{n,i}\ra}_{i=1,2, \cdots, m_n \atop{n\in \NN}}}
\def\alphamn{\cc{\alpha_{n,i}}_{i=1,2, \cdots, m_n \atop{n\in \NN}}}
\def\betamn{\cc{\beta_{n,i}}_{i=1,2, \cdots, m_n \atop{n\in \NN}}}
\def\llim{\lim\limits_{n\to \infty}}
\def\S{\mathfrak{S}}
\def\T{\mathfrak{T}_{\cc{m_n}}}
\def\suml{\sum\limits_{i=1}^{\infty}}

\def\mpsum{\sum\limits_{i=1}^{p}}
\def\mqsum{\sum\limits_{i=1}^{q}}
\def\elii{\ell^2({\scriptscriptstyle{1,2, \cdots, m_n \atop{n\in \NN}}})}
\def\vmn{\mathcal{V}({\scriptscriptstyle{1,2, \cdots, m_n \atop{n\in \NN}}})}

\emph{Frames for Hilbert spaces} were introduced by Duffin and Schaeffer \cite{DS}:

``A sequence $\{x_n\}_{n\in\NN}\subset \nh$ is said to be a \emph{frame}
for a Hilbert space $\nh$ if there exist positive constants $r_1$ and $r_2$
such that
\begin{eqnarray}\label{e1}
r_1\|x\|^2\le \sum\limits_{n=1}^\infty |\langle
x,x_n\rangle|^2\le r_2\|x\|^2, \ \ \text{for all} \ x\in \nh."
\end{eqnarray}
The positive real constants $r_1$ and $r_2$, respectively, are called lower
and upper frame bounds for the frame $\cc{x_n}_{n\in\NN}$.
The inequality (\ref{e1}) is called the \emph{frame inequality} for the frame
$\{x_n\}_{n\in\NN}$. 
A sequence $\cc{x_n}_{n\in\NN}$ in $\nh$ is a \textit{tight frame}    if it is there exist  $r_1,\ r_2$ satisfying inequality $(1)$ with $r_1=r_2$ and a
\emph{Parseval frame} if it is a tight frame with $r_1=r_2=1$. For more details related to frames, one may refer to \cite{C1, CH2}

\medskip

Some important classes  of sequences that are closely related to frames  are the Bessel sequences and Riesz bases. \emph{Bessel sequence}
are the sequences which only satisfies upper inequality of (\ref{e1}). These sequences are in general need not be bases but posses stable reconstruction. On the other-hand, \textit{Riesz basis} being a bounded image of an orthonormal basis in $\nh$ is always a frame for $\nh$.

\medskip
%%%%%%%%%%%%%%%%%%%%%%%%%%%%%%%%%%%%%%%%%%%%%%%%%%%%%%%%%%%%%%%%%%%%%%%%%%%%%%%%%%%%%%%%%%%%%%%%%%%%%%%%%%%%%%%%%%%

After being observed by Daubechies, Grossmann
and Meyer \cite{DGM} in 1986 that, frames provides a stable reconstruction of functions in $L^2({\RR})$, frame theory became popular among researchers. Being redundant in nature, frame representation  has many more benefits  over a basis representation in $L^2(\RR)$, namely signal and image processing \cite{BCE}, filter bank theory \cite{BHF},
wireless communications \cite{HP} and sigma - delta quantization \cite{BPY}. Due to this  prospective, one can  observe frames  as  some kind of  extension  for orthonormal bases in Hilbert spaces.

\medskip

Keeping  more applications in mind, various generalizations of frames have been introduced and studied namely:

\begin{itemize}[leftmargin=0.5cm]
	\item \textit{Fusion frames} by Casazza and Kutyniok  \cite{CK}

\begin{definition}
Let $\{W_i\}_{i \in \NN}$ be a sequence of closed subspaces
in $\nh$ and $\{v_i\}_{i \in \NN}$ be a family of weights, i.e., $v_i > 0$,
for all $i \in \NN$. Then $\{(W_i ,v_i)\}_{i \in \NN}$ is called a \emph{frame
of subspaces} (\emph{fusion frame}), if there exist constants
$0 < C \leq D <\infty$ such that
\begin{eqnarray*}
C \|x\|^2 \leq \sum\limits_{i \in \NN} v_i^2 \|\pi_{W_i}(x)\|^2 \leq
D\|x\|^2, \ \ \text{for all } x \in \nh,
\end{eqnarray*}
where $\pi_{W_i}$ is an orthogonal projection onto the
subspace $W_i$. The constants $C$ and $D$, respectively, are called \emph{lower}
and \emph{upper frame  bounds} for the fusion frame.
\end{definition}

 \item \textit{Pseudo-frames }by Li and Ogawa \cite{LO}

 \begin{definition}[]
Let $\mathcal{X}$ be a closed subspace of a separable Hilbert space $\nh$.
Let $\cc{x_n}_{n\in\NN} $ and  $\cc{x_n^*}_{n\in\NN}$ be    sequences in $\nh$. The sequence $\cc{x_n}_{n\in\NN}$
is a \emph{pseudo frame} for the subspace $\mathcal{X}$  with respect to  $\cc{x_n^*}_{n\in\NN}$
if
\begin{eqnarray*}
 x = \sum\limits_{n=1}^\infty \la x, x_n \ra x_n^*, \ \ \text{for all}\ \ x\in \mathcal{X}.
\end{eqnarray*}
The sequence $\cc{x_n^*}_{n\in\NN}$ is called a \emph{dual pseudo frame}  to $\cc{x_n}_{n\in\NN}$ for the subspace $\mathcal{X}$.
\end{definition}

 \item   \textit{Bounded quasi-projectors} by Fornasier  \cite{F}
\begin{definition}
Let $\nh$ be a Hilbert space and $\mathcal{W}_0$ be a  closed subspace of $\nh$. Let\break $(\mathcal{W}_0, \cc{D_j:\nh \to\nh }_{j\in\NN})$
be a decomposition of $\nh$. Then  a system of
\emph{bounded quasi-projectors} or \emph{a Bessel resolution of the identity} is a set $\mathcal{P}
= \cc{\mathcal{P}_j }_{j\in\NN}$ of operators such that 
\begin{enumerate}[label=({\Roman*}),leftmargin=0.8cm, itemsep=0em]
	\item for each $j\in \NN$, $\mathcal{P}_j :\mathcal{H}\to \mathcal{W}_j(= D_j(\mathcal{W}_0))$, 
\item $\sum\limits_{j\in\NN} \mathcal{P}_j = \mathcal{I}_{\nh}$, in the strong operator topology and
\item there exists positive constant $B$ such that
\begin{eqnarray*}\sum\limits_{j\in\NN} \pp{\mathcal{P}_j(x)}^2  \le B \ \pp{x}^2, \ \ \text{for all} \ x\in \nh.\end{eqnarray*}
\end{enumerate}
\end{definition}
\end{itemize}
The system is called \emph{self-adjoint} and \emph{compatible} with the canonical projections if
\begin{enumerate}[label=({\Roman*}),leftmargin=0.8cm, itemsep=0em]
	\item $\mathcal{P}_j = \mathcal{P}_j^*,$ {for all}\ $j\in\NN$.
	\item $\mathcal{P}_j \circ \pi_{\mathcal{W}_j}=\mathcal{P}_j $, {for all}\ $j\in\NN$.
\end{enumerate}

\noindent

\begin{itemize}[leftmargin=0.5cm]

\item  \textit{$G-$frames }by Sun \cite{S}

 \begin{definition}
Let $\nh$ be a Hilbert space. A sequence
$\{\Lambda_j \in L(\mathcal{H},\mathcal{V}_j): j \in \NN\}$ is said to be a \emph{generalized frame}, or
simply a \emph{$g$-frame} for $\mathcal{H}$ with respect $\{\mathcal{V}_j: j \in \NN\}$
if there exist positive constants $A$ and $B$ such that
\begin{eqnarray*}
A\|x\|^2 \leq \sum\limits_{j \in \NN}\|\Lambda_j x\|^2 \leq B\|x\|^2,  \ \ \text{for all } x \in \mathcal{H},
\end{eqnarray*}
where  $\{\mathcal{V}_j : j \in \NN\}$
is a sequence of subspaces of $\mathcal{H}$.
\end{definition}

\end{itemize}

 Sun also
proved that frames of subspaces (fusion frames), pseudo frames,
bounded quasi-projectors and oblique frames are special cases of
$g-$frames. Fusion frames and  $g-$frames are
also studied in \cite{AR, KA, KM, S}.

\bigskip

Throughout this paper, we will denote $\mathfrak{Q}$ to be a non-commutative field of quaternions,  $\NN$ be the set of natural numbers, $\Q$ be a separable right quaternionic Hilbert space, $\hmn$ be the sequence of separable right quaternionic Hilbert spaces, by the	term ``right linear operator", we mean a ``right $\QQ$-linear operator", $\mathfrak{B}(\Q)$ denotes the set of all bounded (right $\QQ$-linear) operators of $\Q$ and $\mathfrak{B}(\Q, \Qi)$
be the set of all bounded (right $\QQ$-linear) operators from $\Q$ to $\Qi$

The non-commutative field of quaternions $\mathfrak{Q}$ is a four dimensional real algebra with unity. In $\QQ$, $0$ denotes the null element and $1$ denotes the identity with respect to multiplication. It also includes three so-called imaginary units, denoted by $i,j,k$. i.e.,
\begin{eqnarray*}
\mathfrak{Q}=\cc{r_0+r_1i +r_2j +r_3k \ :\ r_0,\ r_1,\ r_2,\  r_3\in \RR}
\end{eqnarray*}
where $i^2=j^2=k^2=-1; \ ij=-ji=k; \ jk=-kj=i$ and $ki=-ik=j$. For each quaternion $q=r_0+r_1i +r_2j +r_3k \in \mathfrak{Q}$, define conjugate of $q$
denoted by $\overline{q}$ as $
\overline{q}=r_0-r_1i -r_2j -r_3k \in \mathfrak{Q}.
$
If $q=r_0+r_1i +r_2j +r_3k$ is a quaternion, then $r_0$ is called the real part of $q$ and $r_1i +r_2j +r_3k$ is called the imaginary part  of $q$. The modulus of a quaternion $q=r_0+r_1i +r_2j +r_3k$ is defined as
\begin{eqnarray*}
|q|=(\overline{q}q)^{1/2} = (q\overline{q})^{1/2}= \sqrt{r_0^2 +r_1^2 +r_2^2 +r_3^2 }.
\end{eqnarray*}
For more literature related to quaternionic Hilbert spaces one may refer to \cite{adler, GMP}

\bigskip

Frames in  separable right quaternionic Hilbert spaces $\Q$ are defined as: 

\begin{definition}[\cite{SS}] Let $\Q$ be a right quaternionic Hilbert space and $\{u_i\}_{i\in \NN}$ be a sequence in $\Q$. Then $\cc{u_i}_{i\in \NN}$ is said to be a \textit{frame} for $V_R(\QQ)$, if there exist two finite real  constants with $0<r_1\le r_2$  such that
	\begin{eqnarray}\label{3.1}
	r_1\|u\|^2\leq\sum_{i\in \NN}|\langle u_i|u\rangle|^2\leq r_2\|u\|^2, \ \text{for all}\ u\in V_R(\QQ).
	\end{eqnarray}
	The positive constants $r_1$ and $r_2$, respectively, are called lower frame and upper frame bounds for the frame $\{u_i\}_{i\in \NN}$. The inequality (\ref{3.1}) is called frame inequality for the frame $\{u_i\}_{i\in I}$. A sequence $\{u_i\}_{i\in \NN}$ is called a \textit{Bessel sequence} for the right quaternionic Hilbert space $\Q$ with bound $r_2$, if $\{u_i\}_{i\in \NN}$ satisfies the right hand side of the inequality (\ref{3.1}). 
	A sequence $\{u_i\}_{i\in \NN}$ is a \textit{tight frame}  for right quaternionic Hilbert space $V_R(\QQ)$ if there exist positive  $r_1$, $r_2$ satisfying inequality (\ref{3.1}) with $r_1=r_2$,  \textit{Parseval frame} if it is tight with $r_1=r_2=1$ and \textit{exact} if it ceases to be a frame in case any one of its element is removed.
\end{definition}
If $\{u_i\}_{i\in \NN}$ is a frame for $V_R(\QQ)$. Then, the right linear operator $T:\ell_2(\QQ)\to V_R(\QQ)$ defined by 
\begin{eqnarray*}T(\{q_i\}_{i\in \NN})=\sum_{i\in \NN}u_iq_i,\ \ \{q_i\}\in\ell_2(\QQ)
\end{eqnarray*}
is called the \textit{(right) synthesis operator} and the adjoint operator $T^*$  is called the \textit{(right) analysis operator} is given by
\begin{eqnarray*}
T^*(u)=\{\langle u_i|u\rangle\}_{i\in \NN},\ u\in V_R(\QQ).
\end{eqnarray*}
Also,  the \textit{(right) frame operator} $S:\Q\rightarrow \Q$ for the frame $\{u_i\}_{i\in I}$ is a right linear operator given by
\begin{eqnarray*}
S(u)=TT^*(u)=T(\{\langle u_i|u\rangle\}_{i\in I})=\sum_{i\in I}u_i\langle u_i|u\rangle,\ u\in \Q.
\end{eqnarray*}

\section*{Various Generalizations of Frames in  Quaternionic Hilbert Spaces}

In this section, we define various generalizations of frames in a right quaternionic Hilbert space. We begin with  the notion of frame of subspaces in a right quaternionic Hilbert space:
\begin{definition}\label{d6}
	Let $\{W_i^R\}_{i \in \NN}$ be a sequence of right closed subspaces of 
	a right separable quaternionic Hilbert space $\Q$ and $\{v_i\}_{i \in \NN}$ be a family of weights, i.e., $v_i > 0$,
	for all $i \in \NN$. Then $\{(W_i^R ,v_i)\}_{i \in \NN}$ is called a \emph{frame
		of subspaces} (\emph{fusion frame}), if there exist constants
	$0 < C \leq D <\infty$ such that
	\begin{eqnarray*}
		C \|x\|^2 \leq \sum\limits_{i \in \NN} v_i^2 \|\pi_{W_i^R}(x)\|^2 \leq
		D\|x\|^2, \ \ \text{for all } x \in \Q,
	\end{eqnarray*}
	where each $\pi_{W_i}$ is an orthogonal projection onto the
	subspace $W_i$. The constants $C$ and $D$, respectively, are called \emph{lower}
	and \emph{upper frame  bounds} for the fusion frame.
\end{definition}

In view of Definition (\ref{d6}) we have the following example:

\begin{example}
	Let $\Q$ be a right separable quaternionic Hilbert space and $\cc{z_i}_{i\in \NN}$ be an orthonormal basis for $\Q$. Define  $W_1^R = [z_1]$, $W_i^R = [z_{i-1}], \ i\ge 2$, and $v_i = 1, \ i\in \NN$. Then  $\{(W_i^R ,v_i)\}_{i \in \NN}$ is a  frame for $\Q$ with lower and upper fusion frame bounds $r_1 = 1$ and $r_2= 2$, respectively.
\end{example}

\medskip

Next, we define the notion of pseudo-frames  using two Bessel sequences in a right quaternionic Hilbert space.

\begin{definition}\label{d8}
	Let $\mathcal{X^R}$ be a right closed subspace of a right separable quaternionic Hilbert space $\Q$.
	Let $\cc{x_i}_{i\in\NN} $ 
	and  $\cc{x_i^*}_{i\in\NN} $ be  sequences in $\Q$. We say that $\cc{x_i}_{i\in\NN}$
	is a \emph{pseudo frame} for the subspace $\mathcal{X^R}$  with respect to  $\cc{x_i^*}_{i\in\NN}$
	if
	\begin{eqnarray*}
		x = \sum\limits_{n=1}^\infty  x_n^*\la  x_n |x \ra, \ \ \text{for all}\ \ x\in \mathcal{X^R}.
	\end{eqnarray*}
	The sequence $\cc{x_n^*}_{n\in\NN}$ is called a \emph{dual pseudo frame}  to $\cc{x_n}_{n\in\NN}$ for the subspace $\mathcal{X^R}$.
\end{definition}

In view of Definition (\ref{d8}) we have a following example:

\begin{example}
	Let $\Q$ be a right separable quaternionic Hilbert space and $\cc{z_i}_{i\in \NN}$ be an orthonormal basis for $\Q$. For each $i \in \NN$, define  ${x_i}= {z_i}$, $x_i^* =  z_{2i-1}$, and $\mathcal{X^R} \subseteq [\cc{x_i^*}_{i\in \NN}]$. Then  $\cc{x_i}_{i\in\NN}$
	is a \emph{pseudo frame} for the subspace $\mathcal{X^R}$  with respect to  $\cc{x_i^*}_{i\in\NN}$.
\end{example}

\medskip
The following definition is an extension of bounded quasi-projectors in right quaternionic Hilbert spaces.

\begin{definition}\label{d10}
Let $\Q$ be a right quaternionic Hilbert space, $\mathcal{W}_0^R$ be a right closed subspace of $\Q$. Let $(\mathcal{W}_0^R, \cc{D_j:\Q\to\Q}_{j\in\NN})$
be a right decomposition of $\Q$. Then a system of 	\emph{bounded quasi-projectors} or \emph{a Bessel resolution of the identity} is a set $\mathcal{P}
	= \cc{\mathcal{P}_j }_{j\in\NN}$ of operators satisfying 
	\begin{enumerate}[label=({\Roman*}),leftmargin=0.8cm, itemsep=0em]
		\item for each $j\in \NN$, $\mathcal{P}_j :\Q\to \mathcal{W}_j(= D_j(\mathcal{W}_0^R))$,
		\item $\sum\limits_{j\in\NN} \mathcal{P}_j = \mathcal{I}_{\Q}$, in the strong operator topology and
		\item there exists a positive constant $r_2$ such that
		\begin{eqnarray*}\sum\limits_{j\in\NN} \pp{\mathcal{P}_j(x)}^2  \le r_2 \ \pp{x}^2, \ \ \text{for all} \ x\in \Q.\end{eqnarray*}
	\end{enumerate}
\end{definition}

\noindent
The system $\cc{\mathcal{P}_j }_{j\in\NN}$ is called \emph{self-adjoint} and \emph{compatible} with the canonical projections if
\begin{enumerate}[label=({\Roman*}),leftmargin=0.8cm, itemsep=0em]
	\item $\mathcal{P}_j = \mathcal{P}_j^*,$ {for all}\ $j\in\NN$.
	\item $\mathcal{P}_j \circ \pi_{\mathcal{W}_j}=\mathcal{P}_j $, {for all}\ $j\in\NN$.
\end{enumerate}

\medskip

In view of Definition (\ref{d10}) we have a following example:

\begin{example}
	Let $\Q$ be a right separable quaternionic Hilbert space and $\cc{z_i}_{i\in \NN}$ be an orthonormal basis for $\Q$. Define  $W_0^R = \Q$, $D_j: \Q \to \Q$ as $D_j\left(x=\sum\limits_{i\in \NN} z_iq_i\right)  = z_jq_j, \ x\in \Q,\ j\in \NN $. Then the system of  operators $\mathcal{P}
	= \cc{\mathcal{P}_j }_{j\in\NN}$  is a	\emph{bounded quasi-projectors} with $r_2=1$. Further, the system of operators $\mathcal{P}
	= \cc{\mathcal{P}_j }_{j\in\NN}$ is  {self-adjoint} and  {compatible} with the canonical projections.
\end{example}

Next, we define   frame of operators in a right quaternionic Hilbert space as follows:
\begin{definition}
Let $\Q$ be a right quaternionic Hilbert space,  $\hmn$ be a sequence of right quaternionic Hilbert spaces, and $\Tmn$ be a sequence of bounded right linear operators. Then $ \tmn$ is 
called a  \emph{frame of operators} for $\Q$ with respect to $\hmn$ if there exist real constants $r_1$ and $r_2$
with $0<r_1\le r_2<\infty$ such that \mathindent6em
\begin{eqnarray}\label{e5}
r_1 \pp{u}^2 \le   \mnsum \pp{\mathfrak{T}_{i}(u)}^2 \le r_2  \pp{u}^2, \ \ \text{for all} \ u\in \Q.
\end{eqnarray}
\end{definition}
The positive constants $r_1$ and $r_2$, respectively, are called \emph{lower} and \emph{upper} bounds
for the frame of operators $\tmn$. A sequence $\Tmn$ is
said to be a \emph{Bessel sequence of
operators} for $\hmn$ with respect to $\Q$ if righthand side of inequality (\ref{e5}) is satisfied.

A frame of operators $\Tmn$ for $\Q$ with respect to $\hmn$  is said to be
\begin{itemize}[leftmargin=.3in]
\item \emph{tight} if it is possible to choose  $r_1,\ r_2$ satisfying inequality (\ref{e5}) with $r_1=r_2$.
\item \emph{Parseval} if it is possible to choose  $r_1,\ r_2$ satisfying inequality (\ref{e5}) with $r_1=r_2=1$.
\end{itemize}

We call $\Tmn$ a  frame of operators for $\Q$
with respect to a right quaternionic Hilbert space $\mathbb{U}^R(\QQ)$ if $\mathbb{H}^R_{i}(\QQ) = \mathbb{U}^R(\QQ), \ i\in \NN$ and simply,  frame of operators  for a right quaternionic Hilbert space $\Q$  if $\mathbb{H}^R_{i}(\QQ) =
 \Q, \ i \in \NN.$
\medskip

Regarding the existence of  frame of operators in a right quaternionic Hilbert
 space, we give the following examples:

\begin{example}\mathindent6em
Let $\Q$ be a right quaternionic Hilbert space and  $N=\cc{z_i}_{i\in \NN}$ be an orthonormal basis for $\Q$.
 Define $\cc{\mathfrak{T}_{i}: \Q \to \QQ}_{i\in\NN}$ by
\begin{eqnarray*}
\mathfrak{T}_{i}(u) = \la z_i|u \ra, \ \ \text{for all} \ u\in \Q, \ \ i \in\NN.
\end{eqnarray*}
 Then $\cc{\mathfrak{T}_{i}: \Q \to \QQ}_{i\in\NN}$ is a  frame of
 operators for $\Q$ with respect $\QQ$.$\hfill\Box$
\end{example}

\begin{example}
Let $\Q$ be a right quaternionic Hilbert space and
 $N=\cc{z_i}_{i\in\NN}$ be a Hilbert basis for $\Q$.
 Define a sequence $\xmn$ in $\Q$ by
\begin{eqnarray*}
\begin{cases}
&u_{1} = z_1, \\
&u_i = z_{i-1}, \ \ i \ge 2,
\end{cases}
\end{eqnarray*}
and $\hmn$ by
\begin{eqnarray*}
\begin{cases}
&\mathbb{H}^R_1(\QQ)=[z_1],\\
&\Qi = [z_{i-1}], \ \ i \ge 2. 
\end{cases}
\end{eqnarray*}
Now, for each $i\in \NN$, define $\mathfrak{T}_{i}: \Q \to \Qi$ by
$
\mathfrak{T}_{i} (u)=  u_{i}\la u_{i}|u \ra,\ \text{for all} \ u\in \Q.
$ 
Then $\Tmn$ is a  frame of operators for $\Q$
with respect to $\hmn$ with  bounds $r_1=1$ and $r_2=2$.$\hfill\Box$
\end{example}

\begin{remark}
Let $\Q$ be a  right quaternionic Hilbert space and $\xmn$ be a frame for $\Q$ such that
$0< \inf\limits_{i\in\NN} \ \pp{u_{i}}\ \le \sup
\limits_{i\in\NN} \ \pp{u_{i}}\ < \infty$. Then 
$\Q$ has a  frame of operators.
\end{remark}

\noindent
Above discussion, give rise to the following observations:\mathindent6em
\begin{enumerate}[leftmargin=0.5cm, label=\textbf{(\Roman*)}]\itemsep4pt
\item  Let $\Q$ be a right quaternionic Hilbert space, $\mathcal{X}^R$ be a right closed subspace of $\Q$ and
$\cc{x_n}_{n\in\NN}$ be a pseudo frame for the subspace  $\mathcal{X}^R$ with respect to $\cc{x_n^*}_{n\in\NN}$.
Then, there exist constants $0 < r_1 \le r_2 < \infty$ such that
\begin{eqnarray*}
r_1\pp{x}^2 \le \suml |\la x_i|x \ra|^2 \le r_2\pp{x}^2,\ \ \text{for all} \ x\in \Q.
\end{eqnarray*}
Define $\cc{\mathfrak{T}_{i}: \mathcal{X}^R \to \QQ}_{i\in\NN}$ by
 $
\mathfrak{T}_{i}(x) = \la  x_{i}|x\ra, \ \text{for all} \ x\in \mathcal{X}^R \ \text{and}\  i  \in\NN.
$
Then, $\cc{\mathfrak{T}_{i}: \mathcal{X} \to \QQ}_{i\in\NN}$ is a frame of operators for $\mathcal{X}^R$ with respect $\QQ$.

\item  Let $\Q$ be a right quaternionic Hilbert space, $\mathcal{W}_0^R$ be a right closed subspace of $\Q$. Let $(\mathcal{W}_0^R, \cc{D_j:\Q\to\Q}_{j\in\NN})$
be a right decomposition of $\Q$ and $\mathcal{P} = \cc{\mathcal{P}_j : \Q \to \mathcal{W}_j^R (= D_j  (\mathcal{W}_0^R))}_{j\in\NN}$ be
a system of  bounded quasi-projectors such that $\mathcal{P}_j =\mathcal{P}_j^*\  \text{and} \
\mathcal{P}_j \circ \pi_{\mathcal{W}_j^R} = \mathcal{P}_j$ for all $j \in \NN$.
 Then there exist
constants $0<r_1\le r_2<\infty$ such that
\begin{eqnarray*}
r_1\pp{x}^2 \ \le\ \sum\limits_{j\in\NN} \ \pp{\mathcal{P}_j(x)}^2\ \le\ r_2\pp{x}^2, \ \ x\in\Q.
\end{eqnarray*}
Define $\cc{\mathfrak{T}_{j}: \Q \to \QQ}_{j\in\NN}$ by $
\mathfrak{T}_{j}(x) = \mathcal{P}_j(x), \ \ \text{for all} \ x\in \Q \ \text{and }  j \in\NN.$ Then $\cc{\mathfrak{T}_{j}: \Q \to \QQ}_{j\in\NN}$ is a frame of
operators for $\Q$ with respect $\QQ$.

\item  Let $\Q$ be a right quaternionic Hilbert space, $\cc{\mathcal{W}_i^R}_{i\in\NN}$ be a sequence of right closed subspace,   	$\cc{{v}_i}_{i\in\NN}$ be a sequence of weights, 
and  $\cc{\mathcal{W}_i^R}_{i\in\NN}$  be a frame of subspaces for $\Q$ with respect 
to $\cc{{v}_i}_{i\in\NN}$. Then
there exist constants $0<r_1\le r_2<\infty$ such that
\begin{eqnarray*}
r_1\pp{x}^2 \ \le \ \suml {v}_i^2 \pp{\pi_{\mathcal{W}_i^R}(x)}^2 \ \le \ r_2\pp{x}^2, \ \ x\in\Q
\end{eqnarray*}
Define $\cc{\mathfrak{T}_{i}: \Q \to \QQ}_{i\in\NN}$ by
 $
\mathfrak{T}_{i}(x) = v_i\ \pi_{\mathcal{W}_i^R}(x), \ \ \text{for all} \ x\in \Q\ \text{and }  i \in\NN.
$ Then $\cc{\mathfrak{T}_{i}: \Q \to \QQ}_{i\in\NN}$ is a frame of
operators for $\Q$ with respect $\CC$.

\end{enumerate}\bigskip

 Let $\hmn$ be a sequence of right quaternionic  Hilbert spaces. Define the space\mathindent6em
\begin{align*}
\mathfrak{H} = \bigoplus_{{i\in\NN}}\Qi = \bigg\{ \xmn : u_{i} \in \Qi \ \text{such that }   \mnsum \ \|u_{i}\|_{\Qi}^2 < \infty  \bigg\}.
\end{align*}
Then $\mathfrak{H}$ is a right quaternionic Hilbert space with the norm given by\mathindent0em
\begin{align*}
\pp{u}_{\mathfrak{H}} \ = \mnsum \ \|u_{i}\|_{\Qi}^2,\ \ \text{for all}\  u \in \mathfrak{H}.
\end{align*}

Next, we give a characterization for a  Bessel sequence of operators.

\begin{theorem}\label{Th27}
Let  $\Q$ be a right quaternionic Hilbert space, $\hmn$ a sequence of 
right quaternionic Hilbert spaces, and $\Tmn$   a sequence of bounded linear operators. Then 
$\tmn$ is a Bessel sequence of operators with Bessel bound $r_2$ for $\Q$ with respect to $\hmn$ if and only if the operator $\S: \Q \to \Q$
defined by
\begin{align*}
\S(x)= \mnsum {\mathfrak{T}_{i}^* \mathfrak{T}_{i}(x)}, \ \ \text{for all}\ x\in\Q
\end{align*}\mathindent4em
is a well-defined bounded linear operator on $\Q$ with $\pp{\S}\le B$.
\end{theorem}

\proof Let $x \in \Q$. Then, for $ p,q \in \NN$ with $  p>q$, we have\mathindent0em
\begin{align*}
\bigg\| \mpsum  \mathfrak{T}_{i}^* \mathfrak{T}_{i}(x) -\mqsum 
\mathfrak{T}_{i}^*  \mathfrak{T}_{i}(x)\bigg\|
&= \sup\limits_{y\in\Q\atop{\pp{y}=1}} \bigg| \bigg\la y \bigg| \mpsum \mathfrak{T}_{i}^* \mathfrak{T}_{i}(x) - \mqsum \mathfrak{T}_{i}^* \mathfrak{T}_{i}(x) \bigg\ra\bigg|\\
%&\le \sup\limits_{y\in\Q\atop{\pp{y}=1}} \bigg[ \mpsum
%| \la y| \mathfrak{T}_{i}^* \mathfrak{T}_{i}(x)  \ra| + \mqsum |\la y| \mathfrak{T}_{i}^* \mathfrak{T}_{i}(x)   \ra|  \bigg]\\
&\le  \sup\limits_{\pp{y}=1}\bigg[ \bigg(\mpsum \|\mathfrak{T}_{i}(x)\|^2\bigg)^{\frac{1}{2}}\bigg(\mpsum \|\mathfrak{T}_{i}(y)\|^2\bigg)^{\frac{1}{2}}\\
&\qquad\qquad  \ +\ \bigg(\mqsum \|\mathfrak{T}_{i}(x)\|^2\bigg)^{\frac{1}{2}}\bigg(\mqsum \|\mathfrak{T}_{i}(y)\|^2\bigg)^{\frac{1}{2}}\bigg]\\
&\le \sqrt{B} \bigg[\bigg(\mpsum \|\mathfrak{T}_{i}(x)\|^2\bigg)^{\frac{1}{2}} +\bigg(\mqsum \|\mathfrak{T}_{i}(x)\|^2\bigg)^{\frac{1}{2}}\bigg].
\end{align*}
Thus $ \mnsum {\mathfrak{T}_{i}^* \mathfrak{T}_{i}(x)} $ exists in $\Q$.
Therefore $\S: \Q \to \Q$ is a well-defined operator on $\Q$.
Further, we compute\mathindent6em
\begin{eqnarray*}
\pp{\S}= \sup\limits_{\pp{x}=1}| \la x| \S(x)\ra|
= \sup_{\pp{x}=1} \mnsum\pp{\mathfrak{T}_{i}(x)}^2
\  \le \ r_2.
\end{eqnarray*}
Conversely, for each $x\in\Q$, we have
\begin{eqnarray*}
 \mnsum\pp{\mathfrak{T}_{i}(x)}^2 =  | \la x| \S(x)\ra|
\le \pp{\S} \ \pp{x}^2 \ \le r_2  \pp{x}^2.
\end{eqnarray*}
Hence, $\Tmn$ is  a  Bessel sequence of operators with Bessel bound $r_2$ for $\Q$ with respect to $\hmn$. \endproof

In view of Theorem \ref{Th27}, if  $\Q$ is a right quaternionic 
Hilbert space, $\hmn$  a sequence of right quaternionic Hilbert spaces, and $\Tmn$  is a frame of operators for $\Q$ with respect to $\hmn$. Then there exist an operator $\S: \Q \to \Q$
defined by
\begin{align*}
\S(x)=  \mnsum {\mathfrak{T}_{i}^* \mathfrak{T}_{i}(x)}, \ \ \text{for all} \ x\in\Q.
\end{align*}\mathindent4em
We call $\S$ to be the \emph{frame operator for the  frame of operators} $\Tmn$.\medskip

Also,  if  $\Q$ is a right quaternionic Hilbert space, $\hmn$   a sequence of right quaternionic Hilbert spaces,
$\Tmn$  a sequence of bounded linear operators, and  $\cc{e_k^{i}}_{k\in\NN}$ is a Hilbert basis in  $\Qi$. Then, for each $i \in  \NN$,  $h^R_i: \Q \to \QQ $ given by $
h^R_i(x) = \la {e_k^{i}}|\mathfrak{T}_{i}(x)\ra,\ \  x\in \Q$ 
 is a right bounded linear functional on $\Q$. Therefore, by Riesz representation Theorem for quaternions, there exists $x_k^{i} \in \Q$ such that
\begin{eqnarray*}
h^R_i(x) = \la x_k^{i}|x\ra , \ \ \text{for all} \ x\in\Q.
\end{eqnarray*}
So, we get
\begin{eqnarray*}
\mathfrak{T}_{i}(x) = \sum\limits_{k=1}^\infty  e_k^{i} \la  e_k^{i}|\mathfrak{T}_{i}(x) \ra\ 
=\sum\limits_{k=1}^\infty e_k^{i} \la x_k^{i}|x\ra\  , \ \ x\in\Q.
\end{eqnarray*}
Thus, for each $x, y \in \Q$, we obtain
\begin{eqnarray*}
\la  \mathfrak{T}_{i}^*(y)|x\ra =  \bigg\la y\bigg| \sum\limits_{k=1}^\infty e_k^{i} \la x_k^{i}|x\ra\ \bigg\ra
= \bigg\la  \sum\limits_{k=1}^\infty x_k^{i} \la e_k^{i} |y\ra  \bigg| x\bigg\ra.
\end{eqnarray*}
Hence
\begin{align*}
 \mathfrak{T}_{i}^*(y)  = \sum\limits_{k=1}^\infty x_k^{i}\la  e_k^{i}|y\ra\   ,\ \ \text{for all} \ y \in \Q.
\end{align*}
 We call the sequence $\cc{x_k^{i}}_{i\in \NN}$  as
 the \emph{sequence induced by the sequence of  right bounded linear operators 
 	 $\Tmn$}.
 \medskip

 In the next result, we show that if
 the sequence induced by $\tmn$ is a   frame for $\Q$, then $\tmn$ is a 
 frame of operators for $\Q$
 and vice-versa.

\begin{theorem}\label{Th28}
Let  $\Q$ be a right quaternionic Hilbert space, $\hmn$  a 
sequence of right quaternionic Hilbert spaces, and $\Tmn$  a sequence of bounded linear operators. Then
$\Tmn$ is a  frame of operators for $\Q$ with respect to $\hmn$ with bounds $r_1$ and $r_2$ if and only if the sequence induced by $\tmn$ is a  frame for $\Q$ with bounds $r_1$ and $r_2$.

Further, the operator $\S$ of  frame of operators $\Tmn$ coincides with the   frame operator
 of the frame induced by $\tmn$.
\end{theorem}

\proof For each $x \in \Q$, we have
\begin{eqnarray*}
 \mnsum\pp{\mathfrak{T}_{i}(x)}^2 =  \mnsum \sum\limits_{k=1}^\infty \ |\la  x_k^{i}|x\ra|^2.
\end{eqnarray*}
 Further, for each $x\in \ \Q$, we compute \mathindent6em
 \begin{eqnarray*}
\S(x)&=&   \mnsum {\mathfrak{T}_{i}^* \mathfrak{T}_{i}(x)}\\
%  &=&  \mnsum\sum\limits_{k=1}^\infty \ x_k^{i}\la e_k^{i}|\mathfrak{T}_{i}(x) \ra  \\
  &=&  \mnsum\sum\limits_{k=1}^\infty \ x_k^{i}\Big\la e_k^{i} \bigg| \sum\limits_{j=1}^\infty e_j^{i} \la x_j^{i}| x \ra\  \Big\ra\\
  &=&  \mnsum\sum\limits_{k=1}^\infty x_k^{i} \la x_k^{i}
  |x \ra\  .
\end{eqnarray*}
Hence, the operator $\S$ of  the frame of operators $\Tmn$ coincides with the  frame operator  of the   frame induced by $\tmn$. $\hfill\Box$\bigskip

Next, we give some properties of the frame operator for frame of operators.

\begin{theorem}
Let  $\Q$ be a right quaternionic Hilbert space,
$\hmn$  a sequence of right quaternionic Hilbert spaces, and $\Tmn$  be a frame of operators for $\Q$ with respect to $\hmn$. Then the frame operator $\S$   for the   frame of operators $\tmn$ is bounded, invertible, self-adjoint and positive.
\end{theorem}

\proof In view of Theorem \ref{Th27}, $\S$ is bounded. For each $x,y \in\Q$, we have\mathindent6em
\begin{eqnarray*}
\la y| \S(x) \ra = \mnsum \la y| {\mathfrak{T}_{i}^* \mathfrak{T}_{i}(x)} \ra
= \mnsum \la  {\mathfrak{T}_{i}^* \mathfrak{T}_{i}(y)}|x  \ra \ = \ \la \S(y) |x \ra
\end{eqnarray*}
Thus, the operator $\S$ is self-adjoint. Also
\begin{eqnarray*}
\big\la x|\S(x) \big\ra =  \mnsum\pp{\mathfrak{T}_{i}(x)}^2, \ \ \text{for all} \ x\in \Q.
\end{eqnarray*}
So, by  inequality (\ref{e5}), we have
\begin{eqnarray*}
r_1  \la x|\mathcal{I}(x) \ra \le \la x|\S(x) \ra \le r_2 \la x|\mathcal{I}(x)\ra,
\ \ \text{for all} \ x \in \Q,
\end{eqnarray*}
where, $\mathcal{I}$ is the identity operator. This gives $r_1 \mathcal{I} \le \S \le r_2 \mathcal{I}.$
Thus $\S$ is a positive operator. Further, since
\begin{align*}
0 \le \mathcal{I}- r_2^{-1}\S \le \frac{r_2-r_1}{r_2}\mathcal{I},
\end{align*}
 $\pp{\mathcal{I}- r_2^{-1}\S}  < 1. $
Hence $\S$ is invertible. \endproof
\mathindent6em

Next, we show that if we have a frame of operators $\Tmn$ with the frame operator $\S$ and bounds $r_1$ and $r_2$, then one can construct another frame of operators with frame operator $\S^{-1}$ and with bounds $r_1^{-1}$  and $r_2^{-1}$.

\begin{theorem}\label{Th30}
Let  $\Q$ be a right quaternionic Hilbert space, $\hmn$  a sequence of right quaternionic Hilbert spaces and $\Tmn$   a  frame of operators
for $\Q$ with respect to $\hmn$ with bounds $r_1$ and $r_2$. Then ${\cc{{\mathfrak{U}_{i}=\mathfrak{T}_{i}\S^{-1}}:
 \Q \to \Qi}_{i\in \NN}}$ is a  frame of operators
for $\Q$ with respect to $\hmn$ with bounds $r_1^{-1}$ and $r_2^{-1}$ and with the frame operator $\S^{-1}$.
\end{theorem}
\proof Since 
\begin{eqnarray*}
x = \S\S^{-1}(x) =  \mnsum {\S^{-1}\mathfrak{T}_{i}^*\mathfrak{T}_{i}(x)},  \ x\in\Q.
\end{eqnarray*}
For each $x \in \Q$, we obtain
\begin{eqnarray*}
x = \mnsum {\mathfrak{T}_{i}^*\mathfrak{U}_{i}(x)}
= \mnsum {\mathfrak{U}_{i}^*\mathfrak{T}_{i}(x)}.
\end{eqnarray*}
Also, we compute \mathindent4em
\begin{eqnarray*}
\mnsum \pp{\mathfrak{U}_{i}(x)}^2 &=& \mnsum \Big\la\mathfrak{T}_{i}\S^{-1}(x)\Big| \mathfrak{T}_{i}\S^{-1}(x)\Big\ra\\
&=& \mnsum \Big\la\S^{-1}(x)\Big|\mathfrak{T}_{ i}^*\mathfrak{T}_{i}\S^{-1}(x)\Big\ra\\
 &\le&  \dfrac{1}{r_1} \pp{x}^2, \ \ \text{for all} \ x\in\Q.
\end{eqnarray*}\mathindent6em
Further, for each $x\in\Q$
\begin{eqnarray*}
\pp{x}^2 &=& \mnsum \Big\la \mathfrak{U}_{i}^*\mathfrak{T}_{i}(x)\Big| x \Big\ra \\
&\le& \bigg( \mnsum \|\mathfrak{T}_{i}(x)\|^2\bigg)^{\frac{1}{2}}\bigg(\mnsum \|\mathfrak{U}_{i}(x)\|^2\bigg)^{\frac{1}{2}}\\
&\le& \sqrt{r_2} \ \pp{x} \ \bigg(\mnsum \|\mathfrak{U}_{i}(x)\|^2\bigg)^{\frac{1}{2}}.
\end{eqnarray*}
This gives
\begin{eqnarray*}
\bigg(\mnsum \|\mathfrak{U}_{i}(x)\|^2\bigg)^{\frac{1}{2}} \ge \dfrac{1}{r_2} \ \pp{x}^2, \ \ \text{for all} \ x\in\Q.
\end{eqnarray*}
Hence $\Umn$ is a  frame of operators with  bounds $1/r_1$  and $1/r_2$. Furthermore, let
$\mathfrak{V}$ denotes the operator for the frame of operators $\ui$. 
Then, for each $x\in\Q$, we have
\begin{eqnarray*}
\S\mathfrak{V}(x) =\mnsum \S\mathfrak{U}_{ i}^*\mathfrak{U}_{ i}(x) 
%&=& \mnsum \S\S^{-1}\mathfrak{T}_{ i}^*\mathfrak{T}_{i}\S^{-1}(x)\\
=  \mnsum \mathfrak{T}_{ i}^*\mathfrak{T}_{i}\S^{-1}(x)\ = \ x.
\end{eqnarray*}
Hence $\mathfrak{V} =\S^{-1}$.\endproof\medskip

In view of Theorem \ref{Th30}, every  frame of operators $\Tmn$ produces a {
 frame of operators}, namely ${\cc{{\mathfrak{U}_{i}=\mathfrak{T}_{i}\S^{-1}}:
 \Q \to \Qi}_{i\in \NN}}$, where $\S$ is the frame operator for the  frame of operators $\Tmn$. We call $\cc{\mathfrak{U}_{i}}_{i\in \NN}$  as \emph{dual of
 frame of operators for $\tmn$}. In such a case $\tmn$ and
$\cc{\mathfrak{U}_{i}}_{i\in \NN}$  are refer to as a pair of dual frame of operators.

\medskip
In the next result, we give a relation between a pair of dual  frame of operators
and their induced sequences.

\begin{theorem}Let 
$\Q$ be a right quaternionic Hilbert space and $\hmn$ be a sequence 
of right quaternionic Hilbert spaces.
If $\Tmn$ and ${\cc{{\mathfrak{U}_{i}}: \Q \to \Qi}}_{i\in \NN}$ forms a pair of dual  frame of operators, then
the sequences induced by them is a pair of dual  frames and vice-versa.
\end{theorem}

\proof Follows from  Theorem \ref{Th28} and Theorem \ref{Th30}. $\hfill\Box$\medskip

Next by using   frame of operators, we construct
 a Parseval  frame of operators.
\begin{theorem}
Let 
$\Q$ be a right quaternionic Hilbert space, $\hmn$  a 
sequence of right quaternionic Hilbert spaces, and $\Tmn$ 
 a  frame of operators
for $\Q$ with respect to $\hmn$. Then  ${\cc{{\mathfrak{R}_{i}=\mathfrak{T}_{i}\S^{-1/2}}:
 \Q \to \Q_{i}}_{i\in \NN}}$ is  a Parseval frame of operators
for $\Q$ with respect to $\hmn$.
\end{theorem}
\proof Straight forward. $\hfill\Box$ \medskip

Finally, in this section we prove that in order to construct a  frame of operator for a right quaternionic Hilbert space, it is enough to construct it for a
dense subset.

\begin{theorem}\label{Th33}
	Let $\Q$ be a right quaternionic Hilbert space, $\hmn$  a sequence of right quaternionic Hilbert spaces and $\Tmn$  a sequence of bounded linear operators  such that there 
	exist positive constants $r_1$ and $r_2$ satisfying
\begin{align*}
r_1 \pp{x}^2\  \le\ \mnsum \|\mathfrak{T}_{i}(x)\|^2 \ \le\ r_2\pp{x}^2, \ \text{for all $x$ in a dense subset $\mathcal{V}$ of $\Q$.}
\end{align*}
Then $\Tmn$ is a 
frame of operators for $\Q$ with respect to $\hmn$ with bounds $r_1$ and $r_2$.
\end{theorem}
\baselineskip14pt
\proof Suppose that there exist $x\in \Q$ such that
$
\mnsum \|\mathfrak{T}_{i}(x)\|^2 > r_2 \pp{x}^2.
$ 
Then, there exists some $n\in \NN$ such that
$
 \mnsum \|\mathfrak{T}_{i}(x)\|^2 > r_2 \pp{x}^2.
$ 
Since $\mathcal{V}$ is dense in $\Q$,  there exist $y\in \mathcal{V}$ such that
$
 \mnsum \|\mathfrak{T}_{i}(y)\|^2> r_2 \pp{y}^2.
$
This a contradiction. Also, note that
\begin{eqnarray}\label{e7}
r_1 \pp{x}^2 \le |\la x|\S(x)\ra|, \ \ \text{for all} \ x\in \mathcal{V}.
\end{eqnarray}
Since $\S$ is bounded and $\mathcal{V}$ is dense in $\Q$,
 we conclude that (\ref{e7})
holds for all $x\in \Q$.$\hfill \Box$

\section*{Stability of  frame of Operators}
In this section, we prove two results related to the stability of  frame of operators
\begin{theorem}
Let 
$\Q$ be a right quaternionic Hilbert space, $\hmn$  a 
sequence of right quaternionic Hilbert spaces,  and
$\Tmn$  a  frame of operators for $\Q$  with respect to $\hmn$  with bounds $r_1$ and $r_2$. Let $\Rmn$
be a sequence of  operators such that there exist constants $\lambda_1,\ \lambda_2,\ \mu \ge 0$ satisfying
\begin{align*}
 \Big(\mnsum \| &(\mathfrak{T}_{i} - \mathfrak{R}_{i})(x)\|^2\Big)^{1/2}\le
 \lambda_1 \Big(\mnsum \| \mathfrak{T}_{i}(x) \|^2\Big)^{1/2} +  \lambda_2 \Big(\mnsum  \|\mathfrak{R}_{i}(x) \|^2\Big)^{1/2}+\ \mu \|{x}\|,\nonumber\\&\qquad\qquad\qquad \qquad\qquad\qquad\qquad\qquad\qquad\qquad\qquad\qquad \ \ \ \ \text{for all  $x \in \Q$}.
\end{align*}
Then $\rmn$ is a frame of operators for $\Q$ with respect to $\hmn$ with bounds
$r_1\ \Big(1- \dfrac{\lambda_1+\lambda_2+ \mu/\sqrt{r_1}}{(1+\lambda_2)}\Big)^2$ and
$r_2\ \Big(1+ \dfrac{\lambda_1+\lambda_2+ \mu/\sqrt{r_1}}{(1-\lambda_2)}\Big)^2$.
\end{theorem}
\proof  For each $x\in\Q$, we have
\begin{align*}
\Big(\mnsum \| &(\mathfrak{T}_{ i} - \mathfrak{R}_{i})(x)\|^2\Big)^{1/2}\le \Bigg(\Big( \lambda_1 + \dfrac{\mu}{\sqrt{r_1}}\Big) \Big(\mnsum \| \mathfrak{T}_{i}(x) \|^2\Big)^{1/2} +\lambda_2 \Big(\mnsum \| \mathfrak{R}_{i}(x) \|^2\Big)^{1/2}\Bigg)
\end{align*}
and
\begin{align*}
 \Big(\mnsum \| (\mathfrak{T}_{i} - \mathfrak{R}_{i})(x)\|^2\Big)^{1/2}\ge
 \Big(\mnsum \| \mathfrak{T}_{i}(x) \|^2\Big)^{1/2} - \Big(\mnsum \| \mathfrak{R}_{i}(x) \|^2\Big)^{1/2}.
\end{align*}
Therefore, we obtain
\begin{align*}
(1+\lambda_2) \Big(\mnsum \| \mathfrak{R}_{i}(x)\|^2\Big)^{1/2}
&\ge\Big(1- \lambda_1- \dfrac{\mu}{\sqrt{r_1}}\Big)  \Big(\mnsum \| \mathfrak{T}_{i}(x) \|^2\Big)^{1/2}\\
&\ge\Big(1- \lambda_1- \dfrac{\mu}{\sqrt{r_1}}\Big)\ \sqrt{r_1}\  \pp{x}, \ \ x\in\Q.
\end{align*}
This gives
\begin{align*}
 \Big(\mnsum \| \mathfrak{R}_{i}(x)\|^2\Big)^{1/2}
\ge\ r_1\ \Big(1- \dfrac{\lambda_1+\lambda_2+ \mu/\sqrt{r_1}}{(1+\lambda_2)}\Big)^2 \ \pp{x}^2.
\end{align*}
On the similar lines, it is easy to verify that\mathindent3em
\begin{eqnarray*}
 \Big(\mnsum \| \mathfrak{R}_{i}(x)\|^2\Big)^{1/2}
\le\ r_2\ \Big(1+ \dfrac{\lambda_1+\lambda_2+ \mu/\sqrt{r_1}}{(1-\lambda_2)}\Big)^2 \ \pp{x}^2, \ x\in\Q.\qquad\ \ \ \ \ \   \qedhere
\end{eqnarray*}\mathindent6em
\medskip

Finally, we discuss  stability of  frame of operators in terms of their conjugate operators.

\begin{theorem}Let 
	$\Q$ be a right quaternionic Hilbert space, $\hmn$ a 
	sequence of right quaternionic Hilbert spaces,  and
$\Tmn$  a frame of operators for $\Q$  with bounds $r_1$ and $r_2$. Assume that $\Rmn$
is a sequence of bounded linear operators such that there exist constants $\lambda,\  \mu \ge 0$ satisfying
$\bigg(\lambda + \dfrac{\mu}{\sqrt{r_1}}\bigg) <1$ and
\begin{align}\label{e8}
\Big\| \mnsum &(\mathfrak{T}_{i}^* - \mathfrak{R}_{i}^*)(x_{i})\Big\|\le\lambda\Bigg\| \mnsum \mathfrak{T}_{i}^*(x_{i}) \Big\| + \mu \mnsum\|{x_{i}}\|^2, \ x_{i} \in \Qi, \ \ i\in \NN.
\end{align}
Then $\Rmn$ is  a frame of operators for $\Q$ with respect to $\hmn$ with bounds
$r_1\ \Big(1- \dfrac{\lambda+ \mu}{\sqrt{r_1}}\Big)^2$ and
$r_2\ \Big(1+ \dfrac{\lambda+ \mu}{\sqrt{r_2}}\Big)^2$.
\end{theorem}

\proof Note that
\begin{eqnarray*}
x_{i} = \sum\limits_{k=1}^\infty e_k^{i}q_k^{i}, \  \cc{q_k^{i}}\in \ell^2(\QQ)\ \text{and } i \in \NN
\end{eqnarray*}
where  $\cc{e_k^{i}}_{k\in\NN}$ is a Hilbert basis in  $\Qi$. For each $k\in \NN$, let $t_k^{i} = \mathfrak{T}_{i}^*(e_k^{i})$ and $r_k^{i} = \mathfrak{R}_{i}^*(e_k^{i})$.
Then, (\ref{e8}) is equivalent to
\begin{eqnarray*}
\bigg\| \mnsum \sum\limits_{k=1}^\infty  (t_k^{i}-r_k^{i})q_k^{i}
 \bigg\| \le \lambda \ \bigg\| \mnsum \sum\limits_{k=1}^\infty t_k^{i} q_k^{i}  \bigg\| + \mu \mnsum \pp{q_k^{i}}^2.
\end{eqnarray*}
Hence the result holds in view of Theorem \ref{Th28} and
 Theorem \ref{Th33}. \endproof

\end{document}